\documentclass[12pt,twoside]{article}

\usepackage{amsmath, amsthm, amssymb}

\date{}

\begin{document}




\centerline{}

\centerline{}

\centerline {\Large{\bf Cantor Primes as Prime-Valued Cyclotomic Polynomials}}

\centerline{}

\newcommand{\mvec}[1]{\mbox{\bfseries\itshape #1}}

\centerline{\bf {Christian Salas}}

\centerline{}

\centerline{The Open University, Walton Hall, Milton Keynes, MK7 6AA, United Kingdom}
\centerline{e-mail: c.p.h.salas@open.ac.uk}

\centerline{}

\newtheorem{Theorem}{\quad Theorem}[section]

\newtheorem{Definition}[Theorem]{\quad Definition}

\newtheorem{Corollary}[Theorem]{\quad Corollary}

\newtheorem{Lemma}[Theorem]{\quad Lemma}

\newtheorem{Example}[Theorem]{\quad Example}

\centerline{}\bigskip

\centerline{\bf {Abstract}}\bigskip

\textit{Cantor primes are primes $p$ such that $1/p$ belongs to the middle-third Cantor set. One way to look at them is as containing the base-3 analogues of the famous Mersenne primes, which encompass all base-2 \emph{repunit} primes, i.e., primes consisting of a contiguous sequence of 1's in base 2 and satisfying an equation of the form $p + 1 = 2^q$. The Cantor primes encompass all base-3 repunit primes satisfying an equation of the form $2p + 1 = 3^q$, and I show that in general all Cantor primes $> 3$ satisfy a closely related equation of the form $2pK + 1 = 3^q$, with the base-3 repunits being the special case $K = 1$. I use this to prove that the Cantor primes $> 3$ are exactly the prime-valued cyclotomic polynomials of the form $\Phi_s(3^{s^j}) \equiv 1$ (mod 4). Significant open problems concern the infinitude of these, making Cantor primes perhaps more interesting than previously realised.}\bigskip

{\bf Keywords:}  \textit{Cantor set, prime numbers, cyclotomic polynomials}

\centerline{}

{\bf Mathematics Subject Classification:}  \emph{11A41, 11R09} 
\section{Introduction}
Any base-$N$ repunit prime $p$ is a cyclotomic polynomial evaluated at N, $\Phi_q(N)$, with $q$ also prime, i.e.,
\begin{equation}
p = \Phi_q(N) = \frac{N^q - 1}{N - 1} = \sum_{k=0}^{q-1} N^k
\end{equation}
It is therefore expressible as a contiguous sequence of 1's in base $N$. For example, $p = 31$ satisfies (1) for $N = 2$ and $q = 5$ and can be expressed as 11111 in base 2. The term \emph{repunit} was coined by A. H. Beiler \cite{BEIL} to indicate that numbers like these consist of repeated units. 
 
The case $N = 2$ corresponds to the famous Mersenne primes on which there is a vast literature \cite{GUY}. They are sequence number A000668 in The Online Encyclopedia of Integer Sequences \cite{SLOANE} and are exactly the prime-valued cyclotomic polynomials of the form $\Phi_s(2) \equiv 3$ (mod 4).

In this note I show that \emph{Cantor primes} can be characterised in a similar way as being exactly the prime-valued cyclotomic polynomials of the form $\Phi_s(3^{s^j}) \equiv 1$ (mod 4). They are primes whose reciprocals belong to the middle-third Cantor set $\mathcal{C}_3$. 

It is easily shown that $\mathcal{C}_3$ contains the reciprocals of all base-3 repunit primes, i.e., those primes $p$ which satisfy an equation of the form $2p + 1 = 3^q$ with $q$ prime. $\mathcal{C}_3$ is a fractal consisting of all the points in $[0, 1]$ which have non-terminating base-3 representations involving only the digits 0 and 2. Rerranging (1) to get the infinite series
\begin{equation}
\frac{1}{p} = \frac{N - 1}{N^q - 1} = \sum_{k=1}^{\infty} \frac{N-1}{N^{qk}}
\end{equation}
and putting $N = 3$ shows that those primes $p$ which satisfy $2p + 1 = 3^q$ are such that $\frac{1}{p}$ can be expressed in base $3$ using only zeros and the digit $2$. This single digit $2$ will appear periodically in the base-$3$ representation of $\frac{1}{p}$ at positions which are multiples of $q$. Since only zeros and the digit $2$ appear in the ternary representation of $\frac{1}{p}$, $\frac{1}{p}$ is never removed in the construction of $\mathcal{C}_3$, so $\frac{1}{p}$ must belong to $\mathcal{C}_3$. 

Base-3 repunit primes are sequence number A076481 in The Online Encyclopedia of Integer Sequences and the exact analogues of the Mersenne primes, i.e., they are the case $N = 3$ in (1). In the next section I show that Cantor primes $> 3$ more generally satisfy a closely related equation of the form $2pK + 1 = 3^q$, with the base-3 repunits being the special case $K = 1$. A subsequent section proves that the Cantor primes $> 3$ are exactly the prime-valued cyclotomic polynomials of the form $\Phi_s(3^{s^j}) \equiv 1$ (mod 4), and a final section considers related open problems.
\section{An Exponential Equation Characterising All Cantor Primes}
\begin{Theorem}
A prime number $p > 3$ is a Cantor prime if and only if it satisfies an equation of the form $2pK + 1 = 3^q$ where $q$ is the order of 3 modulo $p$ and $K$ is a sum of non-negative powers of $3$ each smaller than $3^q$. 
\end{Theorem}

{\it Comment.}
The base-3 repunit primes are then the special case in which $K = 3^0 = 1$. An example is 13, which satisfies $2p + 1 = 3^3$. A counterexample which shows that not all Cantor primes are base-3 repunit primes is 757, which satisfies $26p + 1 = 3^9$ with $K = 3^0 + 3^1 + 3^2 = 13$ and $q = 9$.

{\it Proof.}
Each $x \in \mathcal{C}_3$ can be expressed in ternary form as
\begin{equation}
x = \sum_{k=1}^{\infty} \frac{a_k}{3^k} = 0.a_1a_2\ldots
\end{equation}
where all the $a_k$ are equal to 0 or 2. The construction of $\mathcal{C}_3$ amounts to systematically removing all the points in $[0, 1]$ which cannot be expressed in ternary form with only 0's and 2's, i.e., the removed points all have $a_k = 1$ for one or more $k \in \mathbb{N}$ \cite{OLM}. 

The construction of the Cantor set suggests some simple conditions which a prime number must satisfy in order to be a Cantor prime. If a prime number $p > 3$ is to be a Cantor prime, the first non-zero digit $a_{k_1}$ in the ternary expansion of $\frac{1}{p}$ must be 2. This means that for some $k_1 \in \mathbb{N}$, $p$ must satisfy 
\begin{equation} 
\frac{2}{3^{k_1}} < \frac{1}{p} < \frac{1}{3^{k_1-1}}
\end{equation}
or equivalently
\begin{equation} 
3^{k_1} \in (2p, 3p)
\end{equation}
Prime numbers for which there is no power of 3 in the interval $(2p, 3p)$, e.g., 5, 7, 17, 19, 23, 41, 43, 47, \ldots,\ can therefore be excluded immediately from further consideration. Note that there cannot be any other power of $3$ in the interval (2p, 3p) since $3^{k_1 - 1}$ and $3^{k_1 + 1}$ lie completely to the left and completely to the right of $(2p, 3p)$ respectively. 

If the next non-zero digit after $a_{k_1}$ is to be another 2 rather than a 1, it must be the case for some $k_2 \in \mathbb{N}$ that
\begin{equation}
\frac{2}{3^{k_1 + k_2}} < \frac{1}{p} - \frac{2}{3^{k_1}} < \frac{1}{3^{k_1 + k_2 - 1}} 
\end{equation}
or equivalently
\begin{equation}
3^{k_2} \in \bigg(\frac{2p}{3^{k_1} - 2p}, \frac{3p}{3^{k_1} - 2p}\bigg)
\end{equation}
Thus, any prime numbers for which there is a power of 3 in the interval $(2p, 3p)$ but for which there is no power of 3 in the interval $(\frac{2p}{3^{k_1} - 2p}, \frac{3p}{3^{k_1} - 2p})$ can again be excluded, e.g., 37, 113, 331, 337, 353, 991, 997, 1009. 

Continuing in this way, the condition for the third non-zero digit to be a 2 is
\begin{equation}
3^{k_3} \in \bigg(\frac{2p}{3^{k_2}(3^{k_1} - 2p) - 2p}, \frac{3p}{3^{k_2}(3^{k_1} - 2p) - 2p}\bigg)
\end{equation}
and the condition for the $n$th non-zero digit to be a 2 is
\begin{equation}
3^{k_n} \in \\ \\
\bigg(\frac{2p}{3^{k_{n-1}}(\cdots(3^{k_2}(3^{k_1} - 2p) - 2p)\cdots) - 2p}, \\ \\ \frac{3p}{3^{k_{n-1}}(\cdots(3^{k_2}(3^{k_1} - 2p) - 2p)\cdots) - 2p}\bigg)
\end{equation}
 
The ternary expansions under consideration are all non-terminating, so at first sight it seems as if an endless sequence of tests like these would have to be applied to ensure that $a_k \neq 1$ for any $k \in \mathbb{N}$. However, this is not the case. Let $p$ be a Cantor prime and let $3^{k_1}$ be the smallest power of 3 that exceeds $2p$. Since $p$ is a Cantor prime, both (5) and (9) must be satisfied for all $n$. Multiplying (9) through by $3^{k_1-k_n}$ we get
\begin{equation}
3^{k_1} \in \bigg(\frac{3^{k_1-k_n}\cdot2p}{3^{k_{n-1}}(\cdots(3^{k_2}(3^{k_1} - 2p) - 2p)\cdots) - 2p}, \frac{3^{k_1-k_n}\cdot3p}{3^{k_{n-1}}(\cdots(3^{k_2}(3^{k_1} - 2p) - 2p)\cdots) - 2p}\bigg)
\end{equation}
Since all ternary representations of prime reciprocals $\frac{1}{p}$ for $p > 3$ have a repeating cycle which begins immediately after the point, it must be the case that $k_n = k_1$ for some $n$ in (10). Setting $k_n = k_1$ in (10) we can therefore deduce from the fact that $3^{k_1} \in (2p, 3p)$ and the fact that (10) must be consistent with this for all values of $n$, that all Cantor primes must satisfy an equation of the form
\begin{equation}
3^{k_{n-1}}(\cdots(3^{k_2}(3^{k_1} - 2p) - 2p)\cdots) - 2p = 1
\end{equation}
where $k_1 + k_2 + \cdots + k_{n-1} = q$ is the cycle length in the ternary representation of $\frac{1}{p}$. In other words, $q$ is the order of 3 modulo $p$. By successively considering the cases in which there is only one non-zero term in the repeating cycle, two non-zero terms, three non-zero terms, etc., in (11), and defining
\begin{align*}
d_1 &= q - k_1 \\
d_2 &= q - k_1 - k_2 \\
d_3 &= q - k_1 - k_2 - k_3 \\
\vdots \\
d_n &= q - k_1 - k_2 - \cdots - k_n = 0 
\end{align*}
it is easy to see that (11) can be rearranged as
\begin{equation}
2p\sum_{i=1}^n 3^{d_i} + 1 = 3^q
\end{equation} 
Setting $K = \sum_{i=1}^n 3^{d_i}$, we conclude that every Cantor prime must satisfy an equation of the form $2pK + 1 = 3^q$ as claimed.

Conversely, every prime which satisfies an equation of this form must be a Cantor prime. To see this, note that we can rearrange (12) to get
\begin{equation}
\frac{1}{p} = \frac{2\sum_{i=1}^n 3^{d_i}}{3^q - 1} = 2\sum_{i=1}^n 3^{d_i}\bigg\{\frac{1}{3^q} + \frac{1}{3^{2q}} + \frac{1}{3^{3q}} + \cdots \bigg\}
\end{equation}  
Since $2\sum_{i=1}^n 3^{d_i}$ involves only products of $2$ with powers of $3$ which are each less than $3^q$, (13) is an expression for $\frac{1}{p}$ which corresponds to a ternary representation involving only 2s. Thus, $\frac{1}{p}$ must be in the Cantor set if $2pK + 1 = 3^q$. 
\section{Cantor Primes as Cyclotomic Polynomials}
Let $n$ be a positive integer and let $\zeta_n$ be the complex number $e^{2 \pi i/n}$. The $n^{\text{th}}$ cyclotomic polynomial is defined as
\begin{equation*}
\Phi_n(x) = \prod_{\substack{1 \leq k < n \\
\text{gcd}(k, n) = 1}}
(x - \zeta_n^k)
\end{equation*}
The degree of $\Phi_n(x)$ is $\varphi(n)$ where $\varphi$ is the Euler totient function. There is now a powerful body of theory relating to cyclotomic polynomials and discussions of their basic properties can be found in any textbook on abstract algebra. 

\begin{Lemma}
$x^{(n-1)a} + x^{(n-2)a} + \cdots + x^{2a} + x^a + 1$ is irreducible in $\mathbb{Z}[x]$ if and only if $n = p$ and $a = p^k$ for some prime $p$ and non-negative integer $k$. 
\end{Lemma} 
{\it Proof.} This is proved as Theorem 4 in \cite{GRASSL}.

\begin{Theorem}
A prime number $p > 3$ is a Cantor prime if and only if $p = \Phi_s(3^{s^j}) \equiv 1$ (mod 4) where $s$ is an odd prime and $j$ is a non-negative integer. 
\end{Theorem}

{\it Proof.}
Assume $p$ is a Cantor prime. By Theorem 2.1 we then have 
\begin{equation}
pK = \frac{3^q - 1}{2} = R_q^{(3)}
\end{equation}
where $R_q^{(3)}$ denotes the base-3 repunit consisting of $q$ contiguous units, and $q$ and $K$ are as defined in that theorem. If $q$ is composite, say $q = rs$, we obtain the factorisation
\begin{equation}
R_q^{(3)} = R_r^{(3)}\cdot (3^{(s-1)r} + 3^{(s-2)r} + \cdots + 3^{2r} + 3^r + 1)
\end{equation}
If $q$ is prime we can take $r = 1$. Therefore in both cases at least one factor of $pK$ must be a base-3 repunit. 

If $K = 1$ then $p = R_q^{(3)} = \Phi_s(3)$, since $q$ must be prime in this case. ($R_q^{(3)}$ is composite if $q$ is). If $K > 1$, $p$ is not a base-3 repunit and by Theorem 2.1 K is a sum of powers of 3, so $p$ must be of the general form
\begin{equation}
p = 3^{(s-1)r} + 3^{(s-2)r} + \cdots + 3^{2r} + 3^r + 1
\end{equation} 
for some $s$ and $r$, and $K$ must be a corresponding base-3 repunit $R_r^{(3)}$, otherwise their product could not be $R_{rs}^{(3)}$. But the polynomial in (16) can only be prime if it is irreducible in $\mathbb{Z}[x]$. By Lemma 3.1, this requires $s$ to be a prime number and $r = s^j$ for some non-negative integer $j$, and we therefore have $p = \Phi_s(3^{s^j})$ in this case. We conclude that in all cases we must have $p = \Phi_s(3^{s^j})$ if $p$ is a Cantor prime. Note that $s$ must be an \emph{odd} prime as $\Phi_s(3^{s^j})$ is even for $s = 2$. 

Conversely, suppose that $p = \Phi_s(3^{s^j})$ is a prime number. Then we can multiply it by the base-3 repunit $R_r^{(3)}$ where $r = s^j$ to get the repunit $R_q^{(3)}$ as in (15). Thus, $p$ must satisfy (14) and must therefore be a Cantor prime.

Base-3 repunits are congruent to 0 modulo 4 when they consist of an even number of digits, and to 1 modulo 4 otherwise. Therefore if $p > 3$ is a base-3 repunit prime it must be of the form $4k + 1$. 

If $p$ is prime but not a base-3 repunit, both $r = s^j$ and $q = rs$ in (15) are odd, so both $R_q^{(3)}$ and $R_r^{(3)}$ are base-3 repunits with odd numbers of digits, and thus of the form $4k + 1$. It follows that $p$ is also of the form $4k + 1$ in this case. 

\section{Open Problems}
The infinitude of Cantor primes is currently an open problem shown to be significant in this paper because of the equivalence of Cantor primes and prime-valued cyclotomic polynomials of the form $\Phi_s(3^{s^j})$. 

In the case $j = 0$, it is known that $\Phi_s(3)$ is prime for $s =$ 7, 13, 71, 103, 541, 1091, 1367, 1627, 4177, 9011, 9551, 36913, 43063, 49681, 57917, 483611, and 877843. It seems plausible that there are infinitely many such values of $s$ but this remains to be proved.

The Cantor prime $757 = \Phi_3(3^3)$ is an example with $j > 0$. It is again an open problem to prove there are infinitely many integers $j > 0$ for which $\Phi_s(3^{s^j})$ is prime given a prime $s$, though all such cyclotomic polynomials must be irreducible. 

Previous studies have considered the infinitude of prime-valued cyclotomic polynomials of other types. For example, primes of the form $\Phi_s(1)$ and $\Phi_s(2)$ are studied in \cite{GALLOT}, and other cases are discussed in \cite{DAMIANOU}. 

\centerline{}

{\bf ACKNOWLEDGEMENTS.}

\centerline{} 

I would like to express my gratitude to anonymous reviewers of this manuscript.

\end{document}